\documentclass[10pt,journal,epsfig]{IEEEtran}

\usepackage[dvips]{graphicx}
\usepackage{graphicx}
\usepackage{amssymb}
\usepackage{cite}
\usepackage{subfigure}
\usepackage{amsmath}
\usepackage{algorithm}
\usepackage{algorithmic}
\usepackage{multirow}

\begin{document}

\title{An Iterative Reweighted Method for Tucker Decomposition of Incomplete Multiway Tensors}

\author{Linxiao Yang, Jun Fang, Hongbin Li,~\IEEEmembership{Senior
Member,~IEEE}, and Bing Zeng
\thanks{Linxiao Yang, and Jun Fang are with the National Key Laboratory
of Science and Technology on Communications, University of
Electronic Science and Technology of China, Chengdu 611731, China,
Email: JunFang@uestc.edu.cn}
\thanks{Hongbin Li is
with the Department of Electrical and Computer Engineering,
Stevens Institute of Technology, Hoboken, NJ 07030, USA, E-mail:
Hongbin.Li@stevens.edu}
\thanks{Bing Zeng is with the School of
Electronic Engineering, University of Electronic Science and
Technology of China, Chengdu, 611731, China, Email:
eezeng@uestc.edu.cn}
\thanks{This work was supported in part by the National Science
Foundation of China under Grant 61428103, and the National Science
Foundation under Grant ECCS-1408182.}}

\maketitle



\begin{abstract}
We consider the problem of low-rank decomposition of incomplete
multiway tensors. Since many real-world data lie on an
intrinsically low dimensional subspace, tensor low-rank
decomposition with missing entries has applications in many data
analysis problems such as recommender systems and image
inpainting. In this paper, we focus on Tucker decomposition which
represents an $N$th-order tensor in terms of $N$ factor matrices
and a core tensor via multilinear operations. To exploit the
underlying multilinear low-rank structure in high-dimensional
datasets, we propose a group-based log-sum penalty functional to
place structural sparsity over the core tensor, which leads to a
compact representation with smallest core tensor. The method for
Tucker decomposition is developed by iteratively minimizing a
surrogate function that majorizes the original objective function,
which results in an iterative reweighted process. In addition, to
reduce the computational complexity, an over-relaxed monotone fast
iterative shrinkage-thresholding technique is adapted and embedded
in the iterative reweighted process. The proposed method is able
to determine the model complexity (i.e. multilinear rank) in an
automatic way. Simulation results show that the proposed algorithm
offers competitive performance compared with other existing
algorithms.
\end{abstract}


\section{Introduction}
Multi-dimensional data arise in a variety of applications, such as
recommender systems
\cite{karatzoglou2010multiverse,xiong2010temporal,adomavicius2011context},
multirelational networks
\cite{rai2015leveraging,sutskever2009modelling}, and
brain-computer imaging
\cite{cong2015tensor,cichocki2009nonnegative}. Tensors (i.e.
multiway arrays) provide an effective representation of such data.
Tensor decomposition based on low rank approximation is a powerful
technique to extract useful information from multiway data as many
real-world multiway data are lying on a low dimensional subspace.
Compared with matrix factorization, tensor decomposition can
capture the intrinsic multi-dimensional structure of the multiway
data, which has led to a substantial performance improvement for
harmonic retrieval \cite{haardt2008higher,guo2011variational},
regression/classification
\cite{benetos2008tensor,zhou2014nonnegative,li2014multilinear},
and data completion \cite{adomavicius2011context,rai2014scalable},
etc. Tucker decomposition \cite{tucker1966some} and
CANDECOMP/PARAFAC (CP) decomposition \cite{bro1997parafac} are two
widely used low-rank tensor decompositions. Specifically, CP
decomposes a tensor as a sum of rank-one tensors, whereas Tucker
is a more general decomposition which involves multilinear
operations between a number of factor matrices and a core tensor.
CP decomposition can be viewed as a special case of Tucker
decomposition with a super-diagonal core tensor. It is generally
believed that Tucker decomposition has a better generalization
ability than CP decomposition for different types of data
\cite{cichocki2015tensor}. In many applications, only partial
observations of the tensor may be available. It is therefore
important to develop efficient tensor decomposition methods for
incomplete, sparsely observed data where a significant fraction of
entries is missing.

Low-rank decomposition of incomplete multiway tensors has
attracted a lot of attention over the past few years and a number
of algorithms \cite{liu2009tensor,Giannakis13,mardani2015subspace,
acar2011scalable,chen2014simultaneous,filipovic2014tucker,liu2013tensor,
liu2014generalized,narita2012tensor,zhaobayesian,xu2015bayesian,de2000best}
were proposed via either optimization techniques or probabilistic
model learning. For both CP and Tucker decompositions, the most
challenging task is to determine the model complexity (i.e. the
rank of the tensor) in the presence of missing entries and noise.
It has been shown that determining the CP rank, i.e. the minimum
number of rank-one terms in CP decomposition, is an NP-hard
problem even for a completely observed tensor
\cite{haastad1990tensor}. Unfortunately, many existing methods
require that the rank of the decomposition is specified \emph{a
priori}. To address this issue, a Bayesian method was proposed in
\cite{rai2014scalable} for CP decomposition, where a shrinkage
prior called as the multiplicative gamma process (MGP) was
employed to adaptively learn a concise representation of the
tensor. In \cite{zhaobayesian}, a sparsity-inducing Gaussian
inverse-Gamma prior was placed over multiple latent factors to
achieve automatic rank determination. Besides the above Bayesian
methods, an optimization-based CP decomposition method was
proposed in \cite{Giannakis13,mardani2015subspace}, where the
Frobenius-norm of the factor matrices is used as the rank
regularization to determine an appropriate number of component
tensors.


In addition to the CP rank, another notion of tensor rank is
multilinear rank \cite{kruskal1989rank}, which is defined as the tuple of the ranks of
the mode-$n$ unfoldings of the tensor. Multilinear rank is closely
related to the Tucker decomposition since the multilinear rank is
equivalent to the dimension of the smallest achievable core tensor
in Tucker decomposition \cite{kolda2009tensor}. To search for a low multilinear rank
representation, a tensor nuclear norm, defined as a (weighted)
summation of nuclear norms of mode-$n$ unfoldings, was introduced
to approximate the multilinear rank. Tensor completion and
decomposition can then be accomplished by minimizing the tensor
nuclear norm. Specifically, an alternating direction method of
multipliers (ADMM) was developed in \cite{liu2009tensor,liu2013tensor} to
minimize the tensor nuclear norm with missing data, and
encouraging results were reported on visual data. The success of
\cite{liu2009tensor} has inspired a number of subsequent works
\cite{chen2014simultaneous,filipovic2014tucker,gandy2011tensor,tan2014tensor,liu2014generalized,narita2012tensor}
for tensor completion and decomposition based on tensor nuclear
norm minimization. Nevertheless, the tensor nuclear norm, albeit
effective, is not necessarily the tightest convex envelope of the
multilinear rank \cite{chen2014simultaneous}. Also, the nuclear norm-based methods are
sensitive to outliers and work well only if the tensor is exactly
low multilinear rank.

In this paper, to automatically achieve a concise Tucker
representation, we introduce a notion referred to as the
order-$(N-1)$ sub-tensor and propose a group log-sum penalty
functional to encourage structural sparsity of the core tensor.
Specifically, in the log-sum penalty function, elements in every
sub-tensor of the core tensor along each mode are grouped
together. Minimizing the group log-sum penalty function thus leads
to a structured sparse core tensor with only a few nonzero
order-$(N-1)$ sub-tensors along each mode. By removing the zero
order-$(N-1)$ sub-tensors, the core tensor shrinks and a compact
Tucker decomposition can be obtained. Note that the log-sum
function which behaves like the $\ell_0$-norm is more
sparsity-encouraging than the nuclear norm that is $\ell_1$-norm
applied to the singular values of a matrix. Thus we expect the
group log-sum minimization is more effective than the tensor
nuclear norm-minimization in finding a concise representation of
the tensor. By resorting to a majorization-minimization approach,
we develop an iterative reweighted method via iteratively
decreasing a surrogate function that majorizes the original
log-sum penalty function. The proposed method can determine the
model complexity (i.e. multilinear rank) in an automatic way.
Also, the over-relaxed monotone fast iterative
shrinkage-thresholding technique \cite{yamagishi2011over} is
adapted and embedded in the iterative reweighted process, which
achieves a substantial reduction in computational complexity.

The remainder of this paper is organized as follows. Section
\ref{sec:notation} provides notations and basics on tensors. The
problem of Tucker decomposition with incomplete entries is
formulated as an unconstrained optimization problem in Section
\ref{sec:formulate}. An iterative reweighted method is developed
in Section \ref{sec:iterative} for Tucker decomposition of
incomplete multiway tensors. In Section \ref{sec:accelerate}, the
over-relaxed monotone fast iterative shrinkage-thresholding
technique is adapted and integrated with the proposed iterative
reweighted method, which results in a significant computational
complexity reduction. Simulation results are provided in section
\ref{sec:experiments}, followed by concluding remarks in Section
\ref{sec:conclusion}.

\section{Notations and Basics on Tensors}\label{sec:notation}
We first provide a brief review on tensor decompositions. A tensor
is the generalization of a matrix to higher dimensions, also known
as ways or modes. Vectors and matrices can be viewed as special
cases of tensors with one and two modes, respectively. Throughout
this paper, we use symbols $\otimes$ , $\circ$ and $\ast$ to
denote the Kronecker, outer and Hadamard product, respectively.

Let $\boldsymbol{\mathcal{X}}\in\mathbb{R}^{I_1\times
I_2\times\cdots\times I_N}$ denote an $N$th order tensor with its
$(i_1,\ldots,i_N)$th entry denoted by $\mathcal{X}_{i_1\cdots
i_N}$. Here the order $N$ of a tensor is the number of dimensions.
Fibers are the higher-order analogue of matrix rows and columns.
The mode-$n$ fibers of $\boldsymbol{\mathcal{X}}$ are
$I_n$-dimensional vectors obtained by fixing every index but
$i_n$. Slices are two-dimensional sections of a tensor, defined
by fixing all but two indices.
Unfolding or matricization is an operation that turns a
tensor into a matrix. Specifically, the mode-$n$ unfolding of a
tensor $\boldsymbol{\mathcal{X}}$, denoted as
$\boldsymbol{X}_{(n)}$, arranges the mode-$n$ fibers to be the
columns of the resulting matrix. For notational convenience, we
also use the notation
$\textrm{unfold}_n(\boldsymbol{\mathcal{X}})$ to denote the
unfolding operation along the $n$-th mode. The $n$-mode product of
$\boldsymbol{\mathcal{X}}$ with a matrix
$\boldsymbol{A}\in\mathbb{R}^{J\times I_n}$ is denoted by
$\boldsymbol{\mathcal{X}}\times_n\boldsymbol{A}$ and is of size
$I_1\cdots\times I_{n-1}\times J\times I_{n+1}\times\cdots\times
I_N$, with each mode-$n$ fiber multiplied by the matrix
$\boldsymbol{A}$, i.e.
\begin{align}
\boldsymbol{\mathcal{Y}}=\boldsymbol{\mathcal{X}}\times_n\boldsymbol{A}\Leftrightarrow
\boldsymbol{Y}_{(n)}=\boldsymbol{A}\boldsymbol{X}_{(n)}
\end{align}

\begin{figure}[t]
	\centering
	\includegraphics [width=200pt]{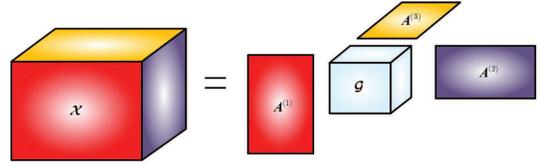}\\
	\caption{The Tucker decomposition of a three-order tensor.}
	\label{fig:tucker}
\end{figure}

The CP decomposition decomposes a tensor into a sum of rank-one
component tensors, i.e.
\begin{align}
\boldsymbol{\mathcal{X}}=
\sum\limits_{r=1}^{R}\lambda_r\boldsymbol{a}_r^{(1)}\circ\boldsymbol{a}_r^{(2)}\circ\cdots\circ\boldsymbol{a}_r^{(N)}
\end{align}
where $\boldsymbol{a}_r^{(n)}\in\mathbb{R}^{I_n}$, `$\circ$'
denotes the vector outer product, and $R$ is referred to as the
rank of the tensor. Elementwise, we have
\begin{align}
\mathcal{X}_{i_1 i_2\cdots i_N}=\sum\limits_{r=1}^{R}\lambda_r
a_{i_1 r}^{(1)}a_{i_2 r}^{(2)}\cdots a_{i_N r}^{(N)}
\end{align}
The Tucker decomposition can be considered as a high order
principle component analysis. It decomposes a tensor into a core
tensor multiplied by a factor matrix along each mode. The Tucker
decomposition of an $N$-th order tensor $\boldsymbol{\mathcal{X}}$
can be written as
\begin{align}
\boldsymbol{\mathcal{X}}&=\boldsymbol{\mathcal{G}}\times_1\boldsymbol{A}^{(1)}
\times_2\boldsymbol{A}^{(2)}\cdots\times_N\boldsymbol{A}^{(N)}\nonumber \\
&=\sum\limits_{r_1=1}^{R_1}\sum\limits_{r_2=2}^{R_2}\cdots\sum\limits_{r_N=1}^{R_N}
\mathcal{G}_{r_1 r_2\cdots
r_N}\boldsymbol{a}^{(1)}_{r_1}\circ\boldsymbol{a}^{(2)}_{r_2}\circ\cdots\circ\boldsymbol{a}^{(N)}_{r_N}\label{tucker}
\end{align}
where $\boldsymbol{\mathcal{G}}\in\mathbb{R}^{R_1\times
R_2\times\cdots\times R_N}$ is the core tensor, and
$\boldsymbol{A}^{(n)}\triangleq
[\boldsymbol{a}_{1}^{(n)}\phantom{0}\ldots\phantom{0}\boldsymbol{a}_{R_n}^{(n)}]\in\mathbb{R}^{I_n\times
R_n}$ denotes the factor matrix along the $n$-th mode (see Fig.\ref{fig:tucker}).

The inner product of two tensors with the same size is defined as
\begin{align}
\langle\boldsymbol{\mathcal{X}},\boldsymbol{\mathcal{Y}}\rangle =
\sum\limits_{i_1=1}^{I_1}\sum\limits_{i_2=1}^{I_2}\cdots
\sum\limits_{i_N=1}^{I_N} x_{i_1 i_2 \dots i_N} y_{i_1 i_2 \dots
i_N}
\end{align}
The Frobenius norm of a tensor $\boldsymbol{\mathcal{X}}$ is
square root of the inner product with itself, i.e.
\begin{align}
\|\boldsymbol{\mathcal{X}}\|_F=\langle\boldsymbol{\mathcal{X}},\boldsymbol{\mathcal{X}}\rangle^{\frac{1}{2}}
\end{align}
Also, for notational convenience, the sequential Kronecker product
of a set of matrices in a reversed order is defined and denoted by
\begin{align}
&\bigotimes_n\boldsymbol{A}^{(n)}=\boldsymbol{A}^{(N)}\otimes\boldsymbol{A}^{(N-1)}\otimes\cdots\otimes\boldsymbol{A}^{(1)}
\nonumber\\
&\bigotimes_{n\ne
k}\boldsymbol{A}^{(n)}=\boldsymbol{A}^{(N)}\otimes\cdots\otimes
\boldsymbol{A}^{(k+1)}\otimes\boldsymbol{A}^{(k-1)}\otimes\cdots\otimes\boldsymbol{A}^{(1)}
\nonumber
\end{align}
An $N$th order tensor $\boldsymbol{\mathcal{X}}$ multiplied by
factor matrices $\{\boldsymbol{A}^{(k)}\}_{k=1}^N$ along each mode
is denoted by
\[
\begin{aligned}
\boldsymbol{\mathcal{X}}\prod\limits_{n=1}^{N}
\times_n\boldsymbol{A}^{(n)}=\boldsymbol{\mathcal{X}}\times_1
\boldsymbol{A}^{(1)}\times_2\boldsymbol{A}^{(2)}\cdots\times_N\boldsymbol{A}^{(N)},
\end{aligned}
\]
while the tensor $\boldsymbol{\mathcal{X}}$ multiplied by the
factor matrices along every mode except the $k$-th mode is denoted
as
\[
\begin{aligned}
&\boldsymbol{\mathcal{X}}\prod\limits_{n\ne k}\times_n\boldsymbol{A}^{(n)}=
\boldsymbol{\mathcal{X}}\times_1\boldsymbol{A}^{(1)}\cdots\times_{k-1}\boldsymbol{A}^{(k-1)}\\
&\qquad\qquad\qquad\qquad\qquad\qquad\times_{k+1}\boldsymbol{A}^{(k+1)}\cdots\times_N\boldsymbol{A}^{(N)}.
\end{aligned}
\]
With these notations, vectorization and unfolding of a tensor
which admits a Tucker decomposition (\ref{tucker}) can be
expressed as
\begin{align}
&\text{vec}(\boldsymbol{\mathcal{X}})=\Big(\bigotimes_n\boldsymbol{A}^{(n)}\Big)
\text{vec}(\boldsymbol{\mathcal{G}})\\
&\text{unfold}_n(\boldsymbol{\mathcal{X}})=\boldsymbol{A}^{(n)}\boldsymbol{G}_{(n)}\Big(\bigotimes_{k\ne
n}\boldsymbol{A}^{(k)}\Big)^T
\end{align}










\section{Problem Formulation}\label{sec:formulate}
Let $\boldsymbol{\mathcal{Y}}\in\mathbb{R}^{I_1\times
I_2\times\cdots\times I_N}$ be an incomplete $N$th order tensor,
with its entry $\mathcal{Y}_{i_1 i_2\dots i_N}$ observed if
$\mathcal{O}_{i_1 i_2\dots i_N}=1$, where
$\boldsymbol{\mathcal{O}}\in\{0,1\}^{I_1\times
I_2\times\cdots\times I_N}$ is a binary tensor of the same size as
$\boldsymbol{\mathcal{Y}}$ and indicates which entries of
$\boldsymbol{\mathcal{Y}}$ are missing or observed. Given the
observed data, our objective is to find a Tucker decomposition
which has a minimum model complexity and meanwhile fits the
observed data, or to be more precise, seek a Tucker representation
such that the data can be represented by a smallest core tensor.
Since the dimension of the smallest achievable core tensor is
unknown \emph{a priori}, we need to develop a method that can
achieve automatic model determination. To this objective, we first
introduce a new notion called as order-$(N-1)$ sub-tensor.

\emph{Definition:} Order-$(N-1)$ sub-tensor is defined as a new
tensor obtained by fixing only one index of the original tensor.
Let $\boldsymbol{\mathcal{Z}}\in\mathbb{R}^{I_1\times
I_2\times\cdots\times I_N}$ be an $N$th order tensor. The $i$th
$(1\le i\le I_n)$ sub-tensor along the $n$th mode of
$\boldsymbol{\mathcal{Z}}$, denoted as
$\boldsymbol{\mathcal{Z}}_{(n,i)}$, is an $(N-1)$th order tensor
of size $I_1\times I_2\times\cdots I_{n-1}\times
I_{n+1}\cdots\times I_N$, and its $(j_1,\ldots,j_{N-1})$th entry
are given by
$\mathcal{X}_{j_1,\ldots,j_{n-1},i,j_{n+1},\ldots,j_{N-1}}$. For
tensors with three modes, i.e. $N=3$, order-$(N-1)$ sub-tensors
reduce to slices, although order-$(N-1)$ sub-tensors are generally
different from slices.


Clearly, $\boldsymbol{\mathcal{Z}}$ consists of $I_n$
order-$(N-1)$ sub-tensors along its $n$th mode. If some
order-$(N-1)$ sub-tensors along the $n$th mode become zero, then
the dimension of $\boldsymbol{\mathcal{Z}}$ along the $n$th mode
is reduced accordingly.
Suppose the data tensor $\boldsymbol{\mathcal{Y}}$ has a Tucker
decomposition
\begin{align}
\boldsymbol{\mathcal{Y}} = \boldsymbol{\mathcal{X}}\prod\limits_{n=1}^{N}\times_n
\boldsymbol{A}^{(n)}
\end{align}
Unfolding $\boldsymbol{\mathcal{Y}}$ along the $n$th mode, we have
\begin{align}
\boldsymbol{Y}_{(n)} =& \boldsymbol{A}^{(n)}\boldsymbol{X}_{(n)}\Big(\bigotimes_{k\ne n}\boldsymbol{A}^{(k)}\Big)^T\nonumber\\
=&\sum\limits_{i=1}^{I_n}\boldsymbol{a}^{(n)}_{\cdot,i}\boldsymbol{x}^{(n)}_{i,\cdot}
\Big(\bigotimes_{k\ne n}\boldsymbol{A}^{(k)}\Big)^T
\end{align}
where $\boldsymbol{a}^{(n)}_{\cdot,i}$ is the $i$th column of $\boldsymbol{A}^{(n)}$
and $\boldsymbol{x}^{(n)}_{i,\cdot}$ is the $i$th row of $\boldsymbol{X}_{(n)}$.
Clearly, $\boldsymbol{x}^{(n)}_{i,\cdot}$ is the vectorization of the
$i$th order $(N-1)$ sub-tensor along the $n$th mode.
If $\boldsymbol{x}^{(n)}_{i,\cdot}$ is a zero vector, i.e. the corresponding
order-$(N-1)$ sub-tensor is a zero tensor, both $\boldsymbol{a}^{(n)}_{\cdot,i}$
and $\boldsymbol{x}^{(n)}_{i,\cdot}$ have no contribution to $\boldsymbol{\mathcal{Y}}$
and can thus be removed.
Inspired by this insight, sparsity can be
enforced over each sub-tensor (along each mode) of the core tensor
such that the observed data can be represented by a structural
sparsest core tensor with only a few nonzero sub-tensors over all
modes.
By removing those zero sub-tensors along each mode (the
associated columns of the factor matrices are disabled and can be
removed as well), the core tensor shrinks to a smaller one and a
compact Tucker decomposition can be obtained. The problem can be
formulated as
\begin{align}
\min_{\boldsymbol{\mathcal{X}},\{\boldsymbol{A}^{(n)}\}}&\quad\sum\limits_{n=1}^{N}\|\boldsymbol{z}_n\|_{0}\nonumber\\
\text{s.t.}&\quad\|\boldsymbol{\mathcal{O}}\ast(\boldsymbol{\mathcal{Y}}-\boldsymbol{\mathcal{X}}
\prod\limits_{n=1}^{N}\times_n\boldsymbol{A}^{(n)})\|_F^2\leq\varepsilon
\label{opt-1}
\end{align}
where $\varepsilon$ is an error tolerance parameter related to
noise statistics, and $\boldsymbol{z}_n$ is an $I_n$-dimensional
vector with its $i$th entry given by
\begin{align}
z_{n,i}\triangleq\|\boldsymbol{\mathcal{X}}_{(n,i)}\|_{F}
\end{align}
It should be noted that since there is usually no knowledge about
the size of smallest core tensor, the dimensions of the core
tensor are predefined to be the same as the original tensor, i.e.
$\boldsymbol{\mathcal{X}}\in\mathbb{R}^{I_1\times
I_2\times\cdots\times I_N}$. The term $\|\boldsymbol{z}_n\|_{0}$
specifies the number of nonzero sub-tensors along the $n$th mode
of tensor $\boldsymbol{\mathcal{X}}$. Since the number of nonzero
sub-tensors along the $n$th mode is equivalent to the dimension of
mode-$n$ of the core tensor, the above optimization yields a
smallest (in terms of sum of dimensions of all modes) core tensor.
The optimization (\ref{opt-1}), however, is an NP-hard problem.
Thus, alternative sparsity-promoting functions which are more
computationally efficient in finding the structural sparse core
tensor are desirable. In this paper, we consider the use of the
log-sum sparsity-encouraging functional. Log-sum penalty function
has been extensively used for sparse signal recovery and was shown
to be more sparsity-encouraging than the $\ell_1$-norm
\cite{candes2008enhancing,chartrand2008iteratively,wipf2010iterative,shen2013exact}.
Replacing the $\ell_0$-norm in (\ref{opt-1}) with the log-sum
functional leads to
\begin{align}
\min_{\boldsymbol{\mathcal{X}},\{\boldsymbol{A}^{(n)}\}}\quad&\sum
\limits_{n=1}^{N}\sum\limits_{i=1}^{I_n}\log(\|\boldsymbol{\mathcal{X}}_{(n,i)}\|_F^2
+\epsilon)\nonumber\\
\text{s.t.}\quad&\|\boldsymbol{\mathcal{O}}\ast(\boldsymbol{\mathcal{Y}}-
\boldsymbol{\mathcal{X}}\prod\limits_{n=1}^{N}\times_n\boldsymbol{A}^{(n)})\|_F^2\leq\varepsilon
\label{opt-2}
\end{align}
where $\epsilon$ is a small positive parameter to ensure the
logarithmic function is well-defined. Note that in our
formulation, coefficients are grouped according to sub-tensors and
different sub-tensors may have overlapping entries. This is
different from the group-LASSO
\cite{tibshirani1996regression,yuan2006model} method in which
entries are grouped into a number of non-overlapping subsets with
sparsity imposed on each subset. Also, to avoid the solution
$\boldsymbol{\mathcal{X}}\to 0$ and
$\{\boldsymbol{A}^{(n)}\}\to\infty$, a Frobenius norm can be
imposed on the factor matrices $\{\boldsymbol{A}^{(n)}\}$. The
optimization (\ref{opt-2}) can eventually be formulated as an
unconstrained optimization problem as follows
\begin{align}
\min_{\boldsymbol{\mathcal{X}},\{\boldsymbol{A}^{(n)}\}}&
L(\boldsymbol{\mathcal{X}},\{\boldsymbol{A}^{(n)}\})=\sum
\limits_{n=1}^{N}\sum\limits_{i=1}^{I_n}\log(\|\boldsymbol{\mathcal{X}}_{(n,i)}\|_F^2
+\epsilon) \nonumber\\
&+\lambda_1\Big\|\boldsymbol{\mathcal{O}}\ast\Big(\boldsymbol{\mathcal{Y}}-\boldsymbol{\mathcal{X}}
\prod\limits_{n=1}^{N}\times_n\boldsymbol{A}^{(n)}\Big)\Big\|_F^2
+\lambda_2\sum\limits_{n=1}^{N}\|\boldsymbol{A}^{(n)}\|_F^2
\label{opt-3}
\end{align}
where $\lambda_1$ is a parameter controlling the tradeoff between
the sparsity of the core tensor and the fitting error, and
$\lambda_2$ is a regularization parameter for factor matrices.
Their choices will be discussed later in our paper.


The above optimization (\ref{opt-3}) can be viewed as searching
for a low multilinear rank representation of the observed data.
Multilinear rank, also referred to as $n$-rank, of an $N$-order
tensor $\boldsymbol{\mathcal{X}}$ is defined as the tuple of the
ranks of the mode-$n$ unfoldings, i.e.
\begin{align}
\text{$n$-rank}\triangleq\{\text{rank}(\boldsymbol{X}_{(1)}),\text{rank}(\boldsymbol{X}_{(2)}),\ldots,
\text{rank}(\boldsymbol{X}_{(N)})\}
\end{align}
It can be shown that $n$-rank is equivalent to the dimensions of
the smallest achievable core tensor in Tucker decomposition\cite{kolda2009tensor}.
Therefore the optimization (\ref{opt-3}) can also be used for
recovery of incomplete low $n$-rank tensors. Existing methods for
low $n$-rank completion employ a tensor nuclear-norm, defined as a
(weighted) summation of nuclear-norms of mode-$n$ unfoldings, to
approximate the $n$-rank and achieve a low $n$-rank
representation. Our formulation, instead, uses the Frobenius-norms
of order-$(N-1)$ sub-tensors to promote a low $n$-rank
representation.











\section{Proposed Iterative Reweighted Method}\label{sec:iterative}
We resort to a bounded optimization approach, also known as the
majorization-minimization (MM) approach \cite{hunter2004tutorial},
to solve the optimization (\ref{opt-3}). The idea of the MM
approach is to iteratively minimize a simple surrogate function
majorizing the given objective function. It can be shown that
through iteratively minimizing the surrogate function, the
iterative process yields a non-increasing objective function value
and eventually converges to a stationary point of the original
objective function.



To obtain an appropriate surrogate function for (\ref{opt-3}), we
first find a suitable surrogate function for the log-sum function.
With the following inequality
\begin{multline}
\log(\|\boldsymbol{\mathcal{X}}_{(n,i)}\|_F^2+\epsilon)
\le\frac{\|\boldsymbol{\mathcal{X}}_{(n,i)}\|_F^2+\epsilon}{\|\boldsymbol{\mathcal{X}}_{(n,i)}^{(t)}\|_F^2+\epsilon}{}\\
+\log(\|\boldsymbol{\mathcal{X}}_{(n,i)}^{(t)}\|_F^2+\epsilon)-1
\end{multline}
we arrive at
$f(\boldsymbol{\mathcal{X}}|\boldsymbol{\mathcal{X}}^{(t)})$
defined in (\ref{Q-function}) is a surrogate function majorizing
the log-sum functional, i.e.
\begin{align}
\sum_{n=1}^{N}\sum_{i=1}^{I_n}&\log(\|\boldsymbol{\mathcal{X}}_{(n,i)}\|_F^2
+\epsilon)\leq
f(\boldsymbol{\mathcal{X}}|\boldsymbol{\mathcal{X}}^{(t)})
\nonumber
\end{align}
with the equality attained when
$\boldsymbol{\mathcal{X}}=\boldsymbol{\mathcal{X}}^{(t)}$, where
\begin{align}
&f(\boldsymbol{\mathcal{X}}|\boldsymbol{\mathcal{X}}^{(t)})
\nonumber\\
\triangleq&\langle\boldsymbol{\mathcal{X}},\boldsymbol{\mathcal{D}}^{(t)}\ast\boldsymbol{\mathcal{X}}\rangle{}
+\sum_{n=1}^{N}\sum_{i=1}^{I_n}\log(\|\boldsymbol{\mathcal{X}}_{(n,i)}^{(t)}\|_F^2+\epsilon)
-\sum_{n=1}^{N}I_n \label{Q-function}
\end{align}
in which $\boldsymbol{\mathcal{D}}^{(t)}$ is a tensor of the same
size of $\boldsymbol{\mathcal{X}}$, with its
$(i_1,i_2,\dots,i_N)$th element given by
\begin{align}
\mathcal{D}_{i_1 i_2\dots
i_N}^{(t)}=\sum_{n=1}^{N}(\|\boldsymbol{\mathcal{X}}_{(n,i_n)}^{(t)}\|_F^2+\epsilon)^{-1}\label{D-update}
\end{align}
Thus we can readily verify that an appropriate surrogate function
majorizing the objective function (\ref{opt-3}) is given as
\begin{multline}
Q(\boldsymbol{\mathcal{X}},\{\boldsymbol{A}^{(n)}\}|\boldsymbol{\mathcal{X}}^{(t)})
=
\lambda_1\Big\|\boldsymbol{\mathcal{O}}\ast\Big(\boldsymbol{\mathcal{Y}}-
\boldsymbol{\mathcal{X}}\prod\limits_{n=1}^{N}\times_n\boldsymbol{A}^{(n)}\Big)\Big\|_F^2{}\\
+\langle\boldsymbol{\mathcal{X}},\boldsymbol{\mathcal{D}}^{(t)}\ast\boldsymbol{\mathcal{X}}\rangle
+\lambda_2\sum\limits_{n=1}^{N}\|\boldsymbol{A}^{(n)}\|_F^2+c
\label{surrogate-function}
\end{multline}
where $c$ is a constant. That is,
\begin{align}
L(\boldsymbol{\mathcal{X}},\{\boldsymbol{A}^{(n)}\}) \leq
Q(\boldsymbol{\mathcal{X}},\{\boldsymbol{A}^{(n)}\}|\boldsymbol{\mathcal{X}}^{(t)})
\end{align}
with the equality attained when
$\boldsymbol{\mathcal{X}}=\boldsymbol{\mathcal{X}}^{(t)}$.

Solving (\ref{opt-3}) now reduces to minimizing the surrogate
function (\ref{surrogate-function}) iteratively. Minimization of
the surrogate function, however, is still difficult since it
involves a joint search over the core tensor
$\boldsymbol{\mathcal{X}}$ and the associated factor matrices
$\{\boldsymbol{A}^{(n)}\}_{n=1}^N$. Nevertheless, we will show
that through iteratively decreasing (not necessarily minimizing)
the surrogate function, the iterative process also results in a
non-increasing objective function value and eventually converges
to a stationary point of
$L(\boldsymbol{\mathcal{X}},\{\boldsymbol{A}^{(n)}\})$. Decreasing
the surrogate function is much easier since we only need to
alternatively minimize the surrogate function
(\ref{surrogate-function}) with respect to each variable while
keeping other variables fixed. Details of this alternating
procedure are provided below.

First, we minimize the surrogate function
(\ref{surrogate-function}) with respect to the core tensor
$\boldsymbol{\mathcal{X}}$, given the factor matrices
$\{\boldsymbol{A}^{(n)}\}$ fixed. The problem reduces to
\begin{align}
\min_{\boldsymbol{\mathcal{X}}}
\lambda_1\Big\|\boldsymbol{\mathcal{O}}\ast\Big(\boldsymbol{\mathcal{Y}}-
\boldsymbol{\mathcal{X}}\prod\limits_{n=1}^{N}\times_n\boldsymbol{A}^{(n)}\Big)\Big\|_F^2{}
+\langle\boldsymbol{\mathcal{X}},\boldsymbol{\mathcal{D}}^{(t)}\ast\boldsymbol{\mathcal{X}}\rangle
\end{align}
Let
$\boldsymbol{x}\triangleq\text{vec}(\boldsymbol{\mathcal{X}})$.
The above optimization can be expressed as
\begin{align}
\min_{\boldsymbol{x}}\lambda_1\Big\|\boldsymbol{\Sigma}\Big(\boldsymbol{y}
-\Big(\bigotimes_{n}\boldsymbol{A}^{(n)}\Big)\boldsymbol{x}\Big)\Big\|_F^2
+\boldsymbol{x}^T\boldsymbol{D}^{(t)}\boldsymbol{x} \label{opt-4}
\end{align}
where
$\boldsymbol{\Sigma}\triangleq\textrm{diag}(\textrm{vec}(\boldsymbol{\mathcal{O}}))$
and
$\boldsymbol{D}^{(t)}\triangleq\textrm{diag}(\textrm{vec}(\boldsymbol{\mathcal{D}}^{(t)}))$.
For notational simplicity, define
\begin{align}
\boldsymbol{H}\triangleq
\Big(\bigotimes_{n}\boldsymbol{A}^{(n)}\Big) \label{H-definition}
\end{align}
The optimal solution to (\ref{opt-4}) can be easily obtained as
\begin{align}
\boldsymbol{x}=(\boldsymbol{H}^T\boldsymbol{\Sigma}\boldsymbol{H}+\lambda_1^{-1}\boldsymbol{D}^{(t)})^{-1}
\boldsymbol{H}^T\boldsymbol{\Sigma}\boldsymbol{y} \label{X-update}
\end{align}

Next, we discuss minimizing the surrogate function
(\ref{surrogate-function}) with respect to the factor matrix
$\boldsymbol{A}^{(n)}$, given that the core tensor
$\boldsymbol{\mathcal{X}}$ and the rest of factor matrices
$\{\boldsymbol{A}^{(k)}\}_{k\neq n}$ fixed. Ignoring terms
independent of $\boldsymbol{A}^{(n)}$ and unfolding the tensor
along the $n$th mode, we arrive at
\begin{align}
\min_{\boldsymbol{A}^{(n)}}\quad&\lambda_1\bigg\|\boldsymbol{O}_{(n)}\ast\bigg(\boldsymbol{Y}_{(n)}-\boldsymbol{A}^{(n)}
\boldsymbol{X}_{(n)}\Big(\bigotimes_{k\neq
n}\boldsymbol{A}^{(k)}\Big)^T\bigg)\bigg\|_F^2 \nonumber\\
&+\lambda_2\|\boldsymbol{A}^{(n)}\|_F^2
\end{align}
Clearly, the optimization can be decomposed into a set of
independent tasks, with each task optimizing each row of
$\boldsymbol{A}^{(n)}$. Specifically, let $\boldsymbol{y}_{n,i}$
denote the $i$th row of $\boldsymbol{Y}_{(n)}$,
$\boldsymbol{a}^{(n)}_i$ denote the $i$th row of
$\boldsymbol{A}^{(n)}$, and
$\boldsymbol{\Sigma}^n_i\triangleq\text{diag}(\boldsymbol{o}_{n,i})$,
with $\boldsymbol{o}_{n,i}$ being the $i$th row of
$\boldsymbol{O}_{(n)}$. The optimization of each row of
$\boldsymbol{A}^{(n)}$ can be written as
\begin{align}
\min_{\boldsymbol{a}^{(n)}_i}\quad&\lambda_1\bigg\|\boldsymbol{\Sigma}^{n}_i\bigg(\boldsymbol{y}_{n,i}-\boldsymbol{a}^{(n)}_i
\boldsymbol{X}_{(n)}\Big(\bigotimes_{k\neq
n}\boldsymbol{A}^{(k)}\Big)^T\bigg)\bigg\|_F^2 \nonumber\\
&+ \lambda_2\|\boldsymbol{a}^{(n)}_i\|_F^2\label{opt-5}
\end{align}
whose optimal solution can be readily given as
\begin{align}
\boldsymbol{a}^{(n)}_i=\lambda_1\boldsymbol{y}_{n,i}\boldsymbol{\Sigma}^n_i
\boldsymbol{\Phi}(\lambda_1\boldsymbol{\Phi}^T\boldsymbol{\Sigma}^n_i
\boldsymbol{\Phi}+\lambda_2\boldsymbol{I})^{-1} \label{A-update}
\end{align}
in which
\[
\begin{aligned}
\boldsymbol{\Phi}\triangleq\Big(\bigotimes_{k\neq n}
\boldsymbol{A}^{(k)}\Big)\boldsymbol{X}_{(n)}^{T}
\end{aligned}
\]
Note that $\boldsymbol{\Phi}$ can be more efficiently calculated
from $\text{unfold}_n(\boldsymbol{\mathcal{X}}\prod_{k\ne
n}\boldsymbol{A}^{(k)})$.

Thus far we have shown how to minimize the surrogate function
(\ref{surrogate-function}) with respect to each variable while
keeping other variables fixed. Given a current estimate of the
core tensor and the associated factor matrices
$\{\boldsymbol{\mathcal{X}}^{(t)},\{(\boldsymbol{A}^{(n)})^{(t)}\}_{n=1}^N\}$,
this alternating procedure is guaranteed to find a new estimate
$\{\boldsymbol{\mathcal{X}}^{(t+1)},\{(\boldsymbol{A}^{(n)})^{(t+1)}\}_{n=1}^N\}$
such that
\begin{align}
Q(\boldsymbol{\mathcal{X}}^{(t+1)},\{(\boldsymbol{A}^{(n)})^{(t+1)}\}|\boldsymbol{\mathcal{X}}^{(t)})
\leq
Q(\boldsymbol{\mathcal{X}}^{(t)},\{(\boldsymbol{A}^{(n)})^{(t)}\}|\boldsymbol{\mathcal{X}}^{(t)})
\label{eqn-1}
\end{align}
In the following, we further show that the new estimate
$\{\boldsymbol{\mathcal{X}}^{(t+1)},\{(\boldsymbol{A}^{(n)})^{(t+1)}\}_{n=1}^N\}$
results in a non-increasing objective function value
\begin{align}
&L(\boldsymbol{\mathcal{X}}^{(t+1)},\{(\boldsymbol{A}^{(n)})^{(t+1)}\})
\nonumber\\
=&
L(\boldsymbol{\mathcal{X}}^{(t+1)},\{(\boldsymbol{A}^{(n)})^{(t+1)}\})
-
Q(\boldsymbol{\mathcal{X}}^{(t+1)},\{(\boldsymbol{A}^{(n)})^{(t+1)}\}|\boldsymbol{\mathcal{X}}^{(t)})
\nonumber\\
&+
Q(\boldsymbol{\mathcal{X}}^{(t+1)},\{(\boldsymbol{A}^{(n)})^{(t+1)}\}|\boldsymbol{\mathcal{X}}^{(t)})
\nonumber\\
\leq&L(\boldsymbol{\mathcal{X}}^{(t)},\{(\boldsymbol{A}^{(n)})^{(t)}\})
-
Q(\boldsymbol{\mathcal{X}}^{(t)},\{(\boldsymbol{A}^{(n)})^{(t)}\}|\boldsymbol{\mathcal{X}}^{(t)})
\nonumber\\
&+
Q(\boldsymbol{\mathcal{X}}^{(t+1)},\{(\boldsymbol{A}^{(n)})^{(t+1)}\}|\boldsymbol{\mathcal{X}}^{(t)})
\nonumber\\
\leq&
L(\boldsymbol{\mathcal{X}}^{(t)},\{(\boldsymbol{A}^{(n)})^{(t)}\})
-
Q(\boldsymbol{\mathcal{X}}^{(t)},\{(\boldsymbol{A}^{(n)})^{(t)}\}|\boldsymbol{\mathcal{X}}^{(t)})
\nonumber\\
&+
Q(\boldsymbol{\mathcal{X}}^{(t)},\{(\boldsymbol{A}^{(n)})^{(t)}\}|\boldsymbol{\mathcal{X}}^{(t)})
\nonumber\\
=&L(\boldsymbol{\mathcal{X}}^{(t)},\{(\boldsymbol{A}^{(n)})^{(t)}\})
\end{align}
where the first inequality follows from the fact that
$Q(\boldsymbol{\mathcal{X}},\{(\boldsymbol{A}^{(n)})\}|\boldsymbol{\mathcal{X}}^{(t)})
-L(\boldsymbol{\mathcal{X}},\{\boldsymbol{A}^{(n)}\})$ attains its
minimum when
$\boldsymbol{\mathcal{X}}=\boldsymbol{\mathcal{X}}^{(t)}$, and the
second inequality comes from (\ref{eqn-1}). We see that through
iteratively decreasing (not necessarily minimizing) the surrogate
function, the objective function
$L(\boldsymbol{\mathcal{X}},\{\boldsymbol{A}^{(n)}\})$ is
guaranteed to be non-increasing at each iteration.

For clarity, we summarize our algorithm as follows.

\begin{center}
\textbf{Iterative Reweighted Algorithm for Tucker Decomposition}
\end{center}

\vspace{0cm} \noindent
\begin{tabular}{lp{7.7cm}}
\hline 1.& Given initial estimates
$\{(\boldsymbol{A}^{(n)})^{(0)}\}$,
$\boldsymbol{\mathcal{X}}^{(0)}$ and
a pre-selected regularization parameter $\lambda_1$ and $\lambda_2$.\\
2.& At iteration $t=0,1,\ldots$: Based on the estimate
$\boldsymbol{\mathcal{X}}^{(t)}$, construct the surrogate function
as depicted in (\ref{surrogate-function}). Search for a new
estimate of $\{(\boldsymbol{A}^{(n)})^{(t+1)}\}$ and
$\boldsymbol{\mathcal{X}}^{(t+1)}$ via (\ref{X-update})
and (\ref{A-update}), respectively. \\
3.& Go to Step 2 if
$\|\boldsymbol{\mathcal{X}}^{(t+1)}-\boldsymbol{\mathcal{X}}^{(t)}\|_F>\varepsilon$,
where $\varepsilon$ is a prescribed tolerance value; otherwise
stop.\\
\hline
\end{tabular}

\vspace{0.3cm}


\emph{Remark:} We discuss the computational complexity of the
proposed method. The main computational task of our proposed
algorithm at each iteration involves calculating a new estimate of
$\boldsymbol{\mathcal{X}}^{(t)}$ and
$\{\boldsymbol{A}^{(n)(t)}\}$. Specifically, the update of the
core tensor $\boldsymbol{\mathcal{X}}$ involves computing an
inverse of a $(\prod_{n}I_n)\times (\prod_{n}I_n)$ matrix (cf.
(\ref{X-update})), which has a computational complexity of order
$O(\prod_{n}I_n^3)$ and scaling polynomially with the data size.
The computational complexity associated with the update of the
$i$th row of $\boldsymbol{A}^{(n)}$ is of order
$O(I_n^3+(\sum_{k\ne n}I_k)\prod_k I_k)$ (cf. (\ref{A-update})),
where the term $O((\sum_{k\ne n}I_k)\prod_k I_k)$ comes from the
computation of $\boldsymbol{\Phi}$ and scales linearly with the
data size, and the term $O(I_n^3)$ is related to the inverse of an
$I_n\times I_n$ matrix. Since all rows of $\boldsymbol{A}^{(n)}$
share a same $\boldsymbol{\Phi}$, the computational complexity of
updating $\boldsymbol{A}^{(n)}$ is of order $O(I_n^4+(\sum_{k\ne
n}I_k)\prod_k I_k)$. We see that the overall computational
complexity at each iteration is dominated by $O(\prod_{n}I_n^3)$,
which scales polynomially with the data size, and makes the
algorithm unsuitable for many real-world applications involving
large dimensions. To address this issue, we resort to a
computationally efficient algorithm, namely, an over-relaxed
monotone fast iterative shrinkage-thresholding algorithm (MFISTA),
to solve the optimization (\ref{opt-4}). A directly calculation of
(\ref{X-update}) is no longer needed and a significant reduction
in computational complexity can be achieved.






\section{A Computationally Efficient Iterative Reweighted Algorithm}\label{sec:accelerate}
It is well known that first order methods based on function values
and gradient evaluations are often practically most feasible
options to solve many large-scale optimization problems. One
famous first order method is the fast iterative
shrinkage-thresholding algorithm (FISTA) \cite{beck2009fast}. It
has a convergence rate of $O(1/k^2)$ for the minimization of the
sum of a smooth and a possibly nonsmooth convex function, where
$k$ denotes the iteration counter. Later on in
\cite{yamagishi2011over}, an over-relaxed monotone FISTA (MFISTA)
was proposed to overcome some limitations inherent in the FISTA.
Specifically, the over-relaxed MFISTA guarantees the monotome
decreasing in the function values, which has been shown to be
helpful in many practical applications. In addition, the
over-relaxed MFISTA admits a variable stepsize in a broader range
than FISTA while keeping the same convergence rate. In the
following, we first provide a brief review of the over-relaxed
MFISTA, and then discuss how to extend the technique to solve our
problem.

\subsection{Review of Over-Relaxed MFISTA}
Consider the general convex optimization problem:
\[
\begin{aligned}
\min_{\boldsymbol{x}}{F(\boldsymbol{x})
    =f(\boldsymbol{x})+g(\boldsymbol{x})}
\end{aligned}
\]
where $f$ is a smooth convex function with the Lipschitz
continuous gradient $L(f)$, and $g$ is a convex but possibly
non-smooth function. The over-relaxed MFISTA scheme is summarized
as follows. Given
$\boldsymbol{x}^{(0)}=\boldsymbol{w}^{(1)},\eta^{(1)}=1$,
$\delta\in(0,2)$ and $\beta\in(0,(2-\delta)/L(f)]$, the sequence
$\{\boldsymbol{x}^{(t)}\}$ is given by
\begin{align}
&\boldsymbol{z}^{(t)}=\text{prox}_{\beta g}(\boldsymbol{w}^{(t)}-
\beta\nabla f(\boldsymbol{w}^{(t)}))\label{MFISTA-1}\\
&\boldsymbol{x}^{(t)}=\arg\min\{F(\boldsymbol{z})|\boldsymbol{z}\in
\{\boldsymbol{z}^{(t)},\boldsymbol{x}^{(t-1)}\}\}\label{comp}\\
&\eta^{(t+1)}=\frac{1+\sqrt{1+4(\eta^{(t)})^2}}{2}\label{MFISTA-3}\\
&\boldsymbol{w}^{(t+1)}=\boldsymbol{x}^{(t)}
+\frac{\eta^{(t)}}{\eta^{(t+1)}}(\boldsymbol{z}^{(t)}-\boldsymbol{x}^{(t)})
+\frac{\eta^{(t)}-1}{\eta^{(t+1)}}(\boldsymbol{x}^{(t)}-\boldsymbol{x}^{(t-1)})\nonumber\\
&\qquad\quad\;+\frac{\eta^{(t)}}{\eta^{(t+1)}}(1-\delta)(\boldsymbol{w}^{(t)}
-\boldsymbol{z}^{(t)})\label{MFISTA-4}
\end{align}
where $\nabla f(\boldsymbol{x})$ denotes the gradient of
$f(\boldsymbol{x})$, and the proximal operator is defined as
\begin{align}
\text{prox}_{\beta g}(\boldsymbol{x}) : \boldsymbol{z} =
\arg\min_{\boldsymbol{z}}\{g(\boldsymbol{z})+\frac{1}{2\beta}\|\boldsymbol{z}
-\boldsymbol{x}\|_2^2\} \label{proximal-operator}
\end{align}
It was proved in \cite{yamagishi2011over} that the sequence
$\{\boldsymbol{x}^{(t)}\}$ is guaranteed to monotonically decrease
the objective function $F(\boldsymbol{x})$ and the convergence
rate is $O(1/k^2)$. Since (\ref{opt-4}) is convex, the
over-relaxed MFISTA can be employed to efficiently solve
(\ref{opt-4}).





\subsection{Solving (\ref{opt-4}) via the Over-Relaxed MFISTA}
Consider the optimization (\ref{opt-4}). Let $f(\boldsymbol{x})$
and $g(\boldsymbol{x})$ respectively represent the data fitting
and regularization terms, i.e.
\begin{align}
f(\boldsymbol{x})&=\lambda_1\Big\|\boldsymbol{\Sigma}\Big(\boldsymbol{y}
-\boldsymbol{H}\boldsymbol{x}\Big)\Big\|_F^2 \nonumber\\
g(\boldsymbol{x})&=\boldsymbol{x}^T\boldsymbol{D}\boldsymbol{x}
\nonumber
\end{align}
Recalling that $\boldsymbol{H}$ is defined in
(\ref{H-definition}). To apply the over-relaxed MFISTA, we need to
compute $\nabla f(\boldsymbol{x})$, $\text{prox}_{\beta
g}(\boldsymbol{x})$, and determine the value of $\beta$. The
gradient of $f(\boldsymbol{x})$ can be easily computed as
\begin{align}
\nabla f(\boldsymbol{x})=2\lambda_1\boldsymbol{H}^T
\boldsymbol{\Sigma}\boldsymbol{H}\boldsymbol{x}
-2\lambda_1\boldsymbol{H}^T\boldsymbol{\Sigma}\boldsymbol{y}
\label{eqn-2}
\end{align}
which can also be expressed in a tensor form as
\begin{align}
\nabla f(\boldsymbol{\mathcal{X}})
=2\lambda_1\bigg(\boldsymbol{\mathcal{O}}\ast\Big(\boldsymbol{\mathcal{X}}\prod
\limits_{n=1}^{N}\times_n\boldsymbol{A}^{(n)}-\boldsymbol{\mathcal{Y}}\Big)\bigg)
\prod\limits_{n=1}^{N}\times_n\boldsymbol{A}^{(n)T}\label{dX}
\end{align}
Such a tensor representation enables a more efficient computation
of $\nabla f(\boldsymbol{x})$. The proximal operation
$\text{prox}_{\beta g}(\boldsymbol{x})$ defined in
(\ref{proximal-operator}) can be obtained as
\begin{align}
\boldsymbol{z}&=\arg\min_{\boldsymbol{z}}\{g(\boldsymbol{z})+
\frac{1}{2\beta}\|\boldsymbol{z}-\boldsymbol{x}\|_2^2\}\nonumber\\
&=\arg\min_{\boldsymbol{z}}\{\boldsymbol{z}^T\boldsymbol{D}
\boldsymbol{z}+\frac{1}{2\beta}\|\boldsymbol{z}-\boldsymbol{x}\|_2^2\}\nonumber\\
&=\arg\min_{\boldsymbol{z}}\{\boldsymbol{z}^T(\boldsymbol{D}+
\frac{1}{2\beta})\boldsymbol{z}-\frac{1}{\beta}\boldsymbol{z}^T\boldsymbol{x}\}\nonumber\\
&=(2\beta\boldsymbol{D}+\boldsymbol{I})^{-1}\boldsymbol{x}\label{prox}
\end{align}
Note that since $\boldsymbol{D}$ is a diagonal matrix, the inverse
of $2\beta\boldsymbol{D}+\boldsymbol{I}$ is easy to calculate. We
now discuss the choice of $\beta$ in the MFISTA. As mentioned
earlier, $\beta$ has to be smaller than $(2-\delta)/L(f)$;
otherwise convergence of the scheme cannot be guaranteed.
Recalling that the Lipschitz continuous gradient $L(f)$ is defined
as any constant which satisfies the following inequality
\begin{align}
\|\nabla f(\boldsymbol{x})-\nabla f(\boldsymbol{y})\|\leq
L(f)\|\boldsymbol{x}-\boldsymbol{y}\| \quad \text{for every
$\boldsymbol{x},\boldsymbol{y}$} \nonumber
\end{align}
where $\|\cdot\|$ denotes the standard Euclidean norm. Hence it is
easy to verify that
\begin{align}
L(f)= 2\lambda_1\lambda_{\text{max}}(\boldsymbol{H}^T
\boldsymbol{\Sigma}\boldsymbol{H})
\end{align}
is a Lipschitz constant of $\nabla f(\boldsymbol{x})$, where
$\lambda_{\text{max}}(\boldsymbol{X})$ denotes the largest
eigenvalue of the matrix $\boldsymbol{X}$. Note that
$\boldsymbol{H}^T \boldsymbol{\Sigma}\boldsymbol{H}$ is of
dimension $(\prod_nI_n)\times (\prod_nI_n)$. Calculation of
$L(f)$, therefore, requires tremendous computational efforts. To
circumvent this issue, we seek an upper bound of $L(f)$ that is
easier to compute. Such an upper bound can be obtained by noticing
that $\boldsymbol{\Sigma}$ is a diagonal matrix with its diagonal
element equal to zero or one
\begin{align}
L(f)=&2\lambda_1\lambda_{\text{max}}(\boldsymbol{H}^T
\boldsymbol{\Sigma}\boldsymbol{H}) \nonumber\\
\stackrel{(a)}{\leq}&2\lambda_1\lambda_{\text{max}}(\boldsymbol{H}^T
\boldsymbol{H}) \nonumber\\
\stackrel{(b)}{=}&2\lambda_1\prod\limits_{n=1}^{N}\lambda_{\text{max}}(\boldsymbol{A}^{(n)T}\boldsymbol{A}^{(n)})
\triangleq\widetilde{L}\label{inequ-1}
\end{align}
where $(a)$ follows from the fact that $\boldsymbol{H}^T
\boldsymbol{H}-\boldsymbol{H}^T \boldsymbol{\Sigma}\boldsymbol{H}$
is positive semi-definite, and $(b)$ comes from the Kronecker
product's properties
\begin{align}
\boldsymbol{H}^T\boldsymbol{H}=\big(\bigotimes_{n}\boldsymbol{A}^{(n)}\big)^T
\big(\bigotimes_{n}\boldsymbol{A}^{(n)}\big)
=\bigotimes_n\big(\boldsymbol{A}^{(n)T}\boldsymbol{A}^{(n)}\big)\nonumber
\end{align}
and
\begin{align}
\text{eig}\Big(\bigotimes_n\big(\boldsymbol{A}^{(n)T}\boldsymbol{A}^{(n)}\big)\Big)
=\bigotimes_n\text{eig}(\boldsymbol{A}^{(n)T}\boldsymbol{A}^{(n)})
\nonumber
\end{align}
in which $\text{eig}(\boldsymbol{X})$ is a vector consisting of
the eigenvalues of matrix $\boldsymbol{X}$. Since
$(2-\delta)/\widetilde{L}\leq (2-\delta)/L(f)$, $\beta$ can be
chosen from $(0,(2-\delta)/\widetilde{L}]$, without affecting the
convergence rate of the over-relaxed MFISTA. The calculation of
$\widetilde{L}$ is much easier than $L(f)$ as the dimension of the
matrix involved in the eigenvalue decomposition has been
significantly reduced.


\emph{Remarks:} We see that the dominant operation in solving
(\ref{opt-4}) via the over-relaxed MFISTA is the evaluation of
gradient (\ref{dX}), which has a computational complexity of order
$O((\sum_{n}I_n)\prod_{n}I_n)$ that scales linearly with the data
size. Thus a significant reduction in computational complexity is
achieved as compared with a direct calculation of
(\ref{X-update}). In addition, our proposed iterative reweighted
method only needs to decrease (not necessarily minimize) the
surrogate function at each iteration. Therefore when applying the
over-relaxed MFISTA to solve (\ref{opt-4}), there is no need to
wait until convergence is achieved. Only a few iterations are
enough since the over-relaxed MFISTA guarantees a monotome
decreasing in the function values. This further reduces the
computational complexity of the proposed algorithm.



For clarity, we now summarize the proposed computationally
efficient iterative reweighted algorithm as follows.

%
%
%

\begin{algorithm}
    \renewcommand{\algorithmicrequire}{\textbf{Input:}}
    \renewcommand\algorithmicensure {\textbf{Output:}}
    \caption{Iterative Re-weighted Algorithm For Incomplete Tensor Decomposition}
    \begin{algorithmic}[1]
        \REQUIRE $\boldsymbol{\mathcal{Y}}$, $\boldsymbol{\mathcal{O}}$,
        $\delta$, $\lambda_1$ and $\lambda_2$
        \ENSURE $\boldsymbol{\mathcal{X}}$, $\{\boldsymbol{A}^{(n)}\}_{n=1}^N$, $\boldsymbol{\mathcal{Y}}$ and
        multilinear rank
        \STATE Initialize $\boldsymbol{\mathcal{X}}$,
        $\{\boldsymbol{A}^{(n)}\}$,
        $\boldsymbol{\mathcal{D}}$.
        \WHILE {not converge}
        \STATE Calculate $\widetilde{L}$ using (\ref{inequ-1}) and select $\beta$ from $(0,(2-\delta)/\widetilde{L}]$
        \STATE Set $\boldsymbol{x}^{(0)}=\text{vec}(\boldsymbol{\mathcal{X}})$
        \FOR{$t=1$ to $t_{\text{max}}$}
        \STATE Calculate the gradient of $f(\text{tensor}(\boldsymbol{x}^{(t)}))$ using (\ref{dX}) and the proximal operation
        $\text{prox}_{\beta g}(\boldsymbol{x})$ using (\ref{prox})
        \STATE Update $\boldsymbol{x}^{(t)}$ using (\ref{MFISTA-1}), (\ref{comp}), (\ref{MFISTA-3}), (\ref{MFISTA-4})
        \ENDFOR
        \STATE Set $\boldsymbol{\mathcal{X}}=\text{tensor}(\boldsymbol{x}^{(t_{\text{max}})})$
        \FOR{$n=1$ to $N$}
        \FOR{$i=1$ to $I_n$}
        \STATE Update $\boldsymbol{a}_i^{(n)}$ using (\ref{A-update})
        \ENDFOR
        \ENDFOR
        \STATE Remove the zero order-$(N-1)$ sub-tensors of $\boldsymbol{\mathcal{X}}$ and corresponding
        columns of $\{\boldsymbol{A}^{(n)}\}$.
        \ENDWHILE
        \STATE Reconstruct $\boldsymbol{\mathcal{Y}}$ using estimated $\boldsymbol{\mathcal{X}}$ and $\{\boldsymbol{A}^{(n)}\}$
        \STATE Estimate multilinear rank by count the nonzero order-$(N-1)$ sub-tensors of estimated
        $\boldsymbol{\mathcal{X}}$ along each mode
    \end{algorithmic}
    \label{algorithm}
\end{algorithm}

\begin{table*}[t]
    \centering
    \caption{Comparison of different algorithms for tensor completion (synthetic and Chemometrics Data)}
    \begin{tabular}{|c|c|c|c|c|c|c|c|c|c|}
        \hline
        &  & \multicolumn{2}{|c|}{Synthetic (CP)} & \multicolumn{2}{|c|}{Synthetic (Tucker)} &
        \multicolumn{2}{|c|}{Amino} & \multicolumn{2}{|c|}{Flow} \\
        \hline
        &  & $50\%$ & $80\%$ & $50\%$ & $80\%$ & $50\%$ & $80\%$ & $50\%$ & $80\%$ \\
        \hline
        \multirow{4}{*}{LRTI}
        & NMSE & 0.0676 & 0.1384 & 0.0723 & 0.1425 & 0.0603 & 0.0955 & 0.1068 & 0.1119 \\
        \cline{2-10}
        & Rank & 6 & 6 & 7 & 7 & 3 & 3 & 7 & 6 \\
        \cline{2-10}
        & Std(R) & 0 & 0 & 0.3162 & 0.4830 & 0.5270 & 0 & 1.7512 & 0.5164 \\
        \cline{2-10}
        & runtime & 4.0049 & 1.5719 & 22.1448 & 19.5644 & 10.6141 & 9.3226 & 43.9058 & 32.7794 \\
        \hline
        \multirow{2}{*}{HaLRTC}
        & NMSE & 0.3141 & 0.6757 & 0.2746 & 0.4786 & 0.2513 & 0.2840 & 0.2459 & 0.2639 \\
        \cline{2-10}
        & runtime & 10.3533 & 6.7849 & 8.1829 & 11.6388 & 59.4361 & 63.5765 & 20.7997 & 23.9988 \\
        \hline
        \multirow{2}{*}{WTucker}
        & NMSE & 0.1623 & 0.3747 & 0.1049 & 0.1949 & 0.1310 & 0.2544 & 0.0895 & 0.1625 \\
        \cline{2-10}
        & runtime & 88.5976 & 202.1250 & 51.3877 & 86.8674 & 42.1004 & 64.0108 & 33.1926 & 66.2170 \\
        \hline
        \multirow{4}{*}{IRTD}
        & NMSE & 0.0660 & 0.1157 & 0.0500 & 0.0857 & 0.0580 & 0.0880 & 0.0809 & 0.0862 \\
        \cline{2-10}
        & $n$-Rank & $(6,6,6)$ & $(6,6,6)$ & $(3,4,5)$ & $(3,4,5)$ & $(3,3,3)$ & $(3,3,3)$ & $(4,5,3)$ & $(4,6,4)$ \\
        \cline{2-10}
        & Std(R) & $(0,0,0)$ & $(0,0,0)$ & $(0,0,0)$ & $(0,0,0)$ & $(0,0,0)$ & $(0,0,0)$ & $(0,0,0)$ & $(0,0.4830,0.5164)$ \\
        \cline{2-10}
        & runtime & 20.8098 & 39.0132 & 32.3008 & 51.8319 & 30.5949 & 121.7123 &
        46.0298 & 157.4853 \\
        \hline
    \end{tabular}
    \label{table1}
\end{table*}


\section{Simulation Results}\label{sec:experiments}
In this section, we conduct experiments to illustrate the
performance of our proposed iterative reweight Tucker
decomposition method (referred to as IRTD). In our simulations, we
set $\delta=0.1$, $\beta=(2-\delta)/\widetilde{L}(f)$ and
$\lambda_2=1$. In fact, our proposed algorithm is insensitive to
the choices of these parameters. The choice of $\lambda_1$ is more
critical than the others, and a suitable choice of $\lambda_1$
depends on the noise level and the data missing ratio. Empirical
results suggest that stable recovery performance can be achieved
when $\lambda_1$ is set in the range $[0.5, 2]$. The factor
matrices and core tensor are initialized by decomposing the
observed tensor (the missing elements are set to zero) with high
order singular value decomposition \cite{de2000multilinear}. In
our algorithm, the over-relaxed MFISTA performs only two
iterations to update the core tensor, i.e. $t_{\text{max}}=2$. We
compare our method with several existing state-of-the-art tensor
decomposition/completion methods, namely, a CP decomposition-based
tensor completion method (also referred to as the low rank tensor
imputation (LRTI)) which uses the Frobenius-norm of the factor
matrices as the rank regularization \cite{Giannakis13}, a tensor
nuclear-norm based tensor completion method \cite{liu2013tensor}
which is also referred to as the high accuracy low rank tensor
completion (HaLRTC) method, and a Tucker factorization method
based on pre-specified multilinear rank \cite{filipovic2014tucker}
which is referred to as the WTucker method. It should be noted
that the LRTI requires to set a regularization parameters
$\lambda$ to control the tradeoff between the rank and the data
fitting error, the HaLRTC method is unable to provide an explicit
multilinear rank estimate, and the WTucker method requires an
over-estimated multilinear rank. All the parameters used for
competing algorithms are tuned carefully to ensure the best
performance is achieved.


\begin{figure*}[t]
    \centering
    \includegraphics [width=100pt]{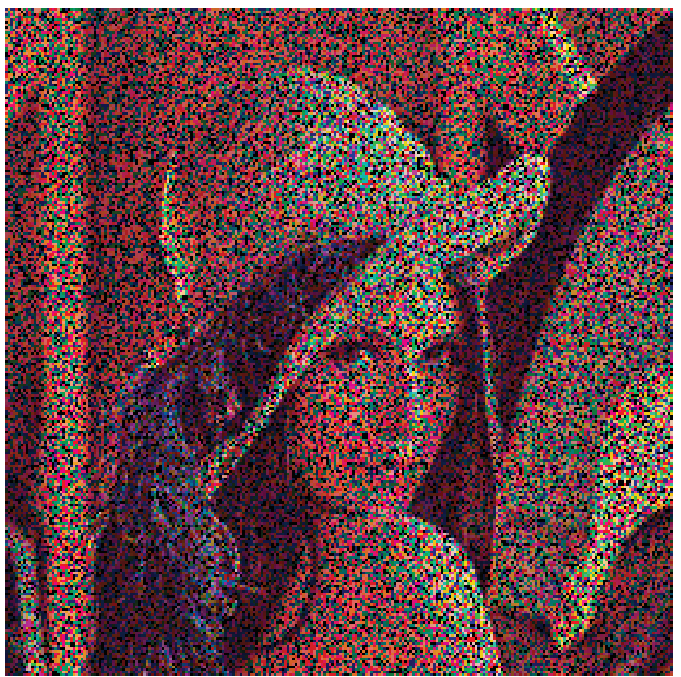}
    \includegraphics [width=100pt]{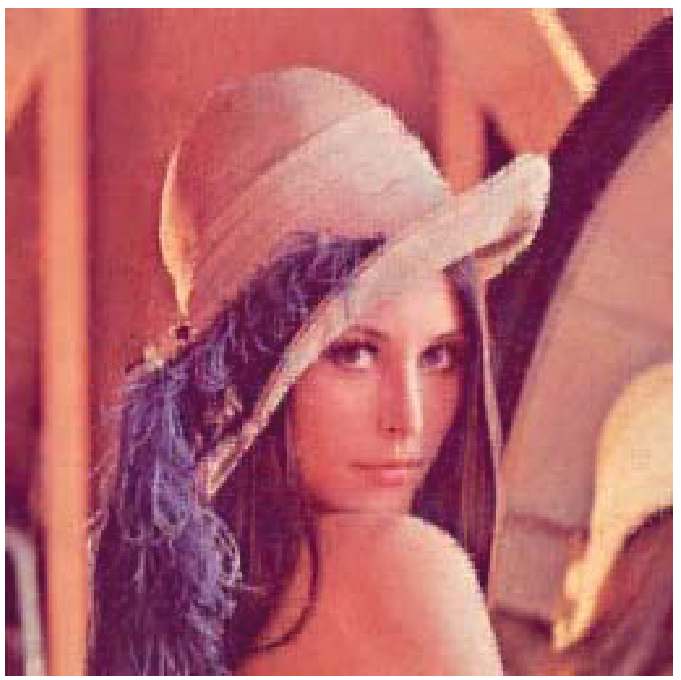}
    \includegraphics [width=100pt]{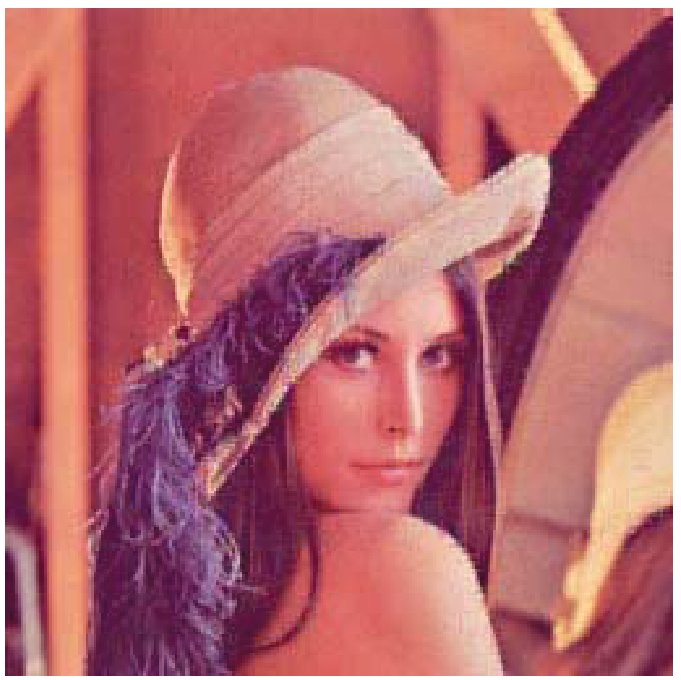}
    \includegraphics [width=100pt]{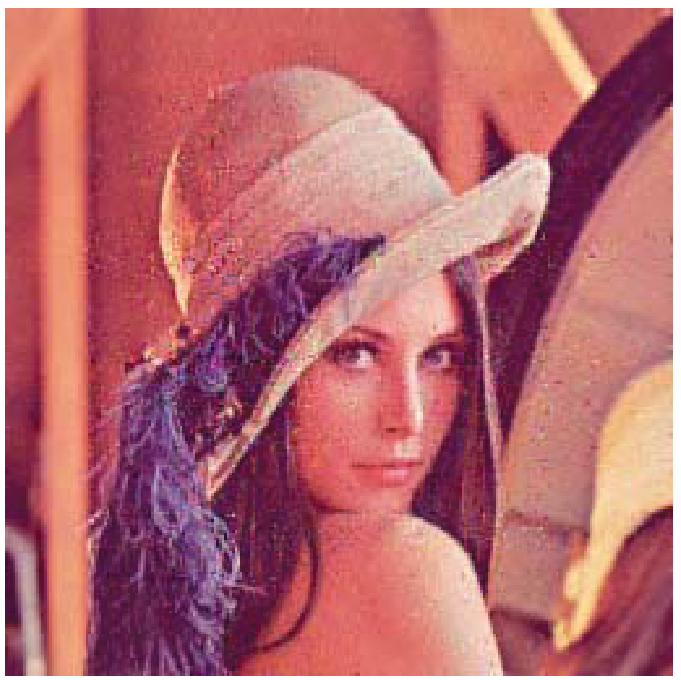}
    \includegraphics [width=100pt]{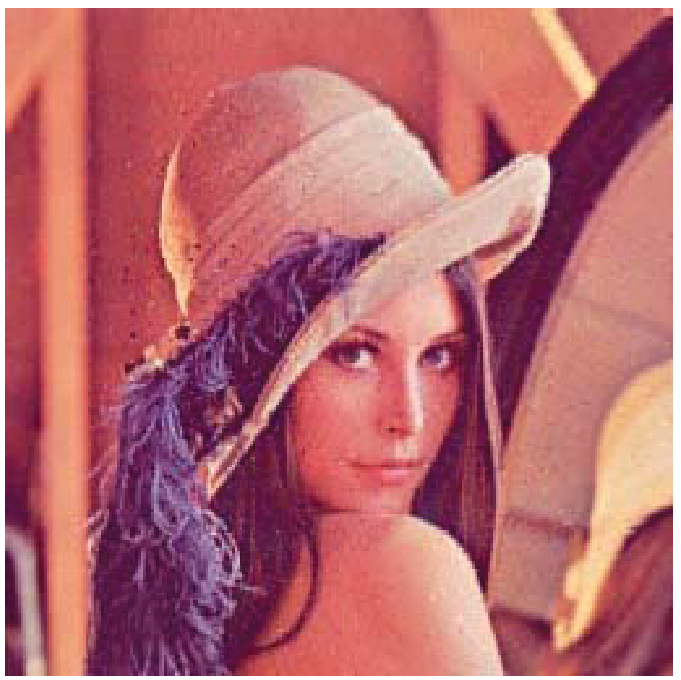}\\

    \includegraphics [width=100pt]{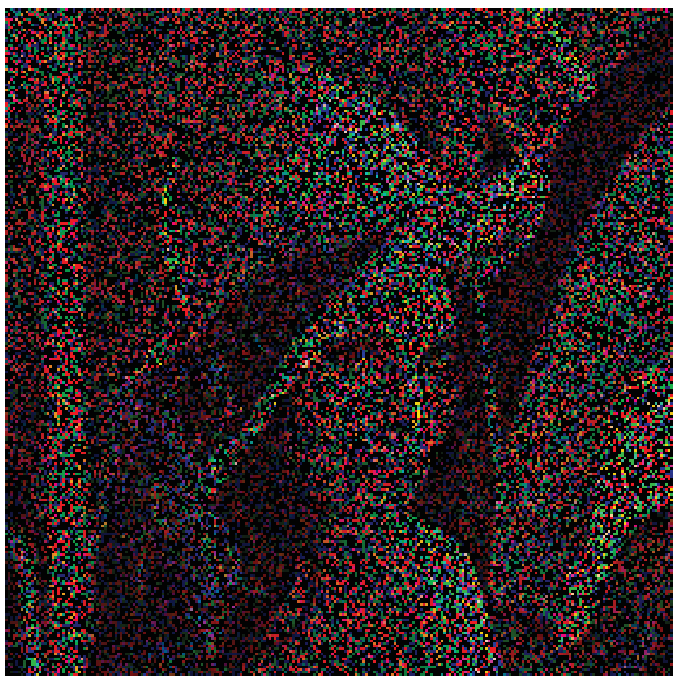}
    \includegraphics [width=100pt]{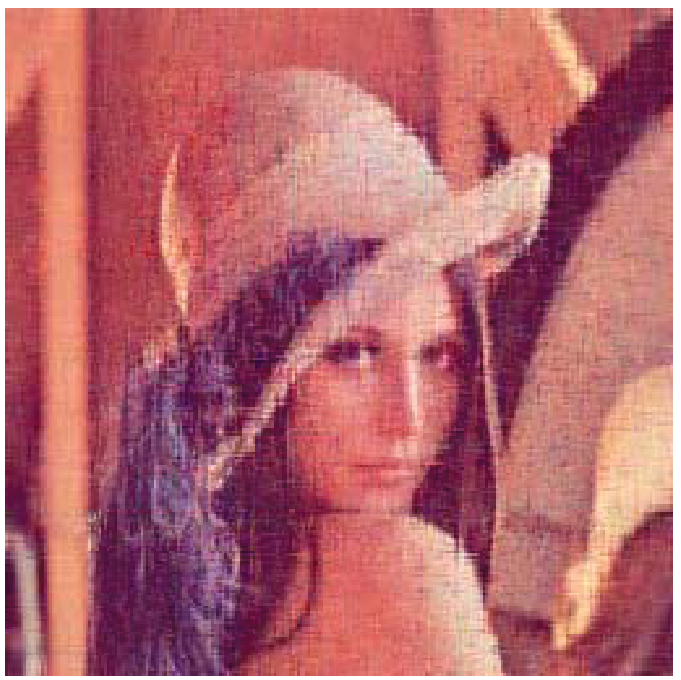}
    \includegraphics [width=100pt]{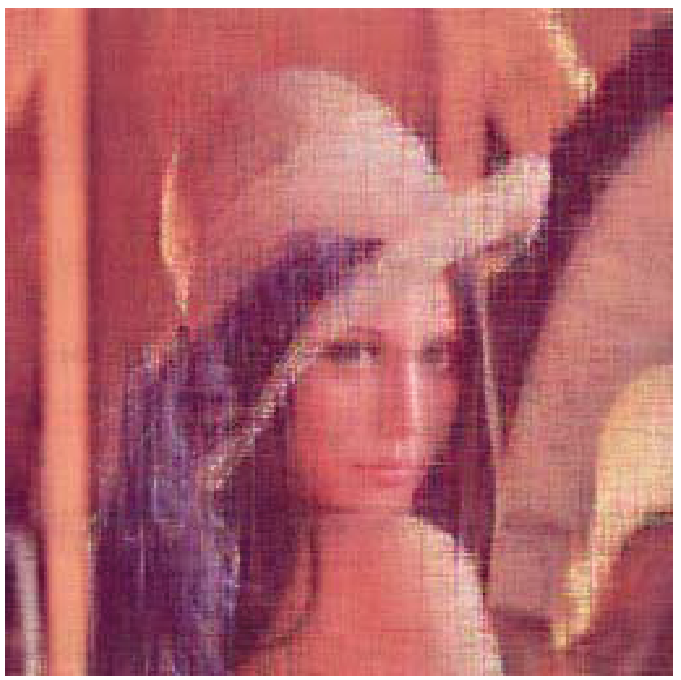}
    \includegraphics [width=100pt]{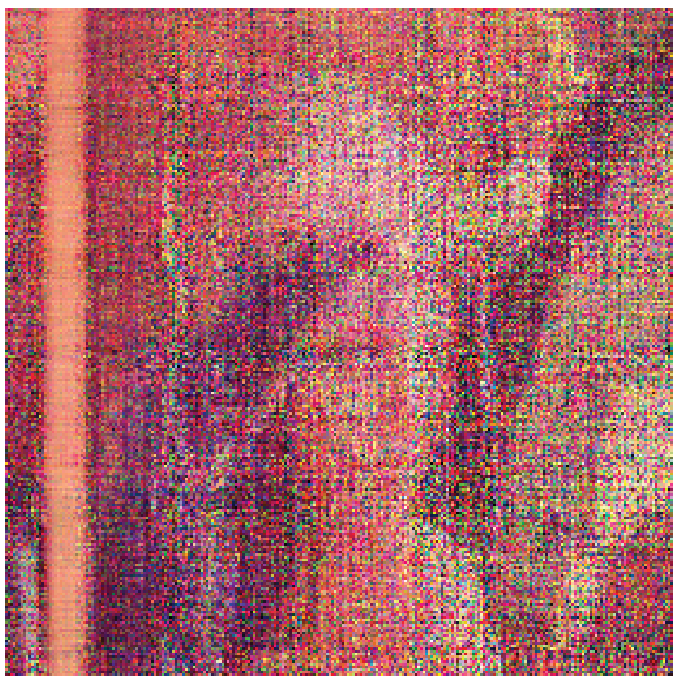}
    \includegraphics [width=100pt]{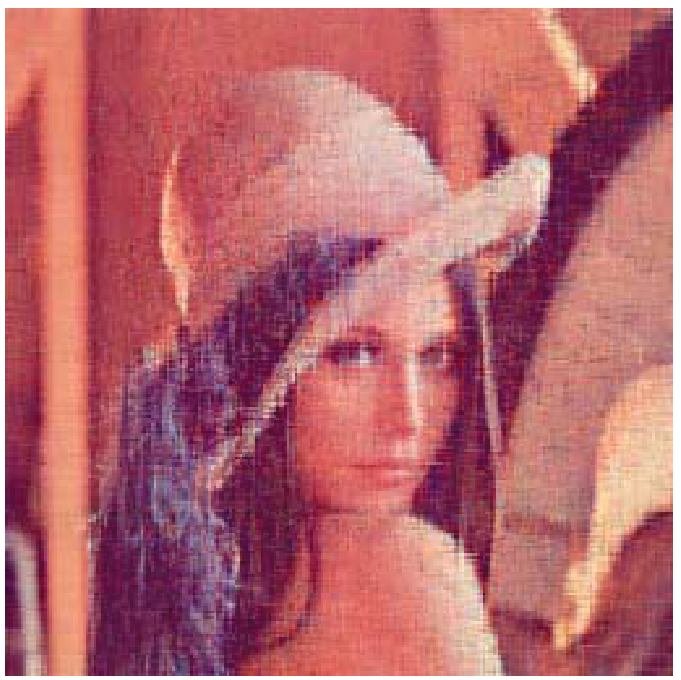}\\

    \includegraphics [width=100pt]{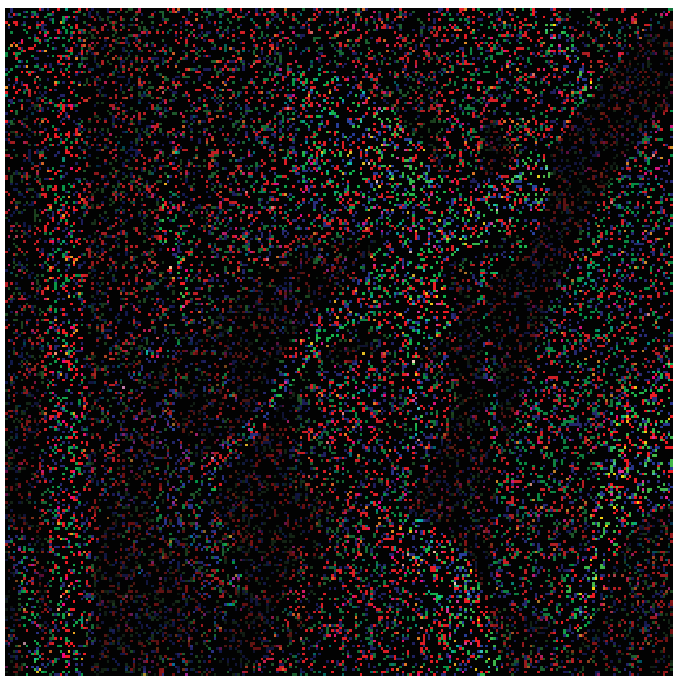}
    \includegraphics [width=100pt]{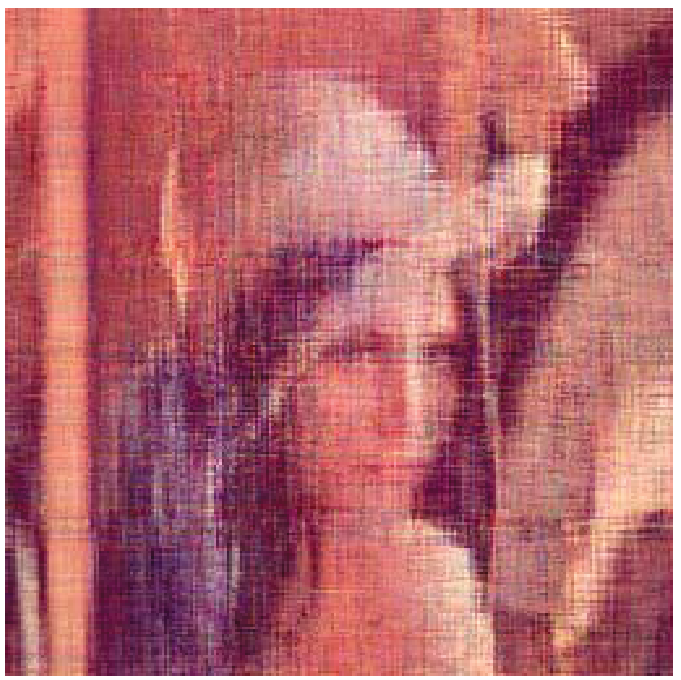}
    \includegraphics [width=100pt]{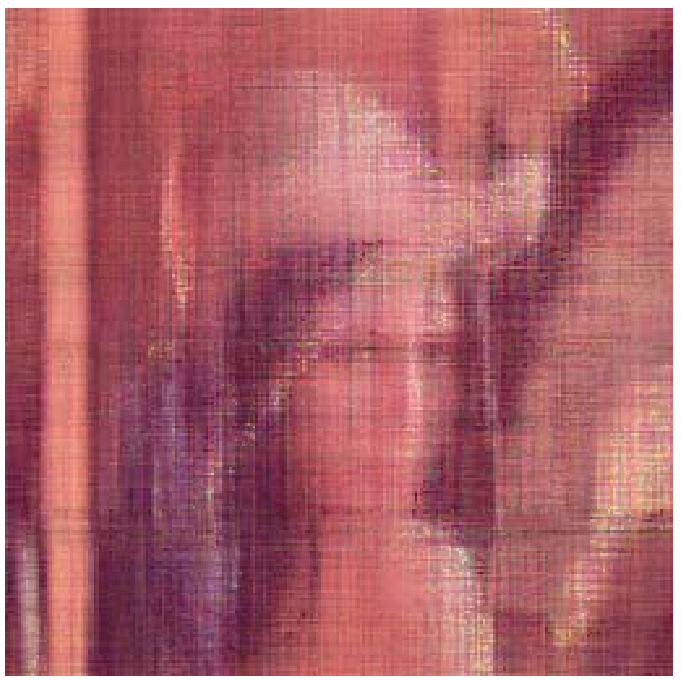}
    \includegraphics [width=100pt]{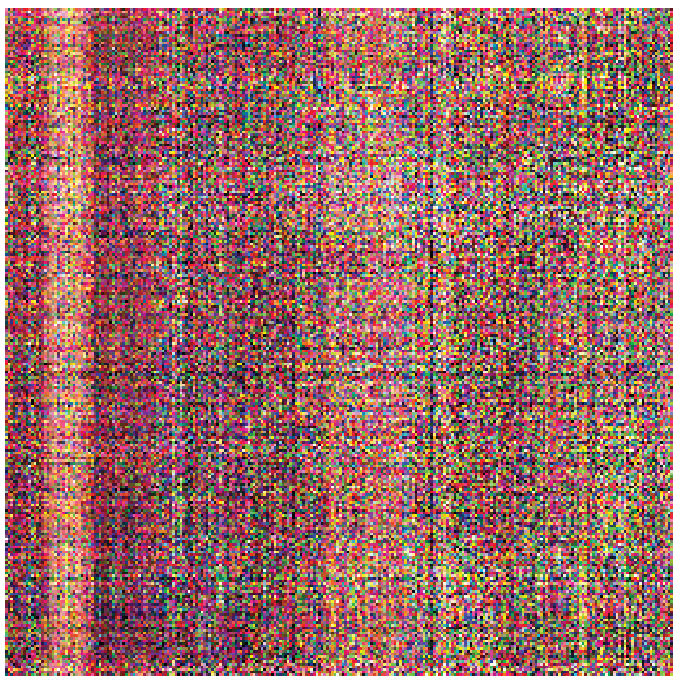}
    \includegraphics [width=100pt]{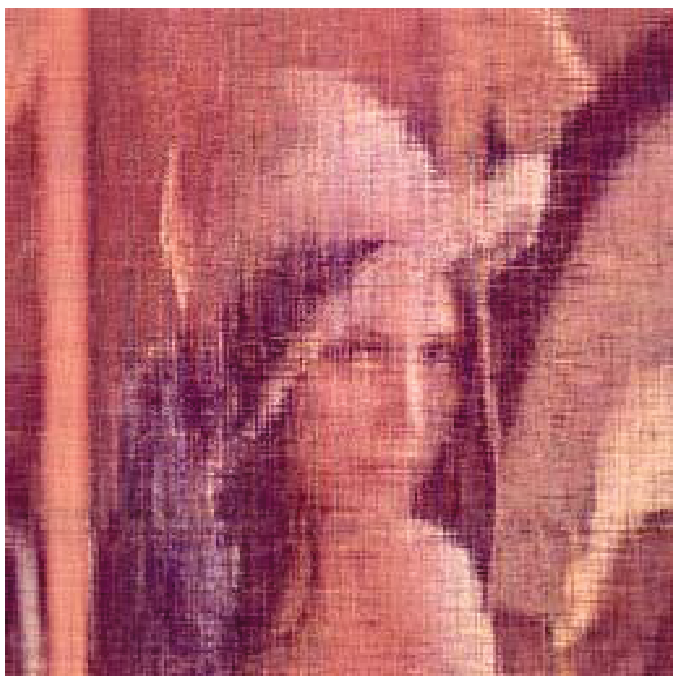}\\

    \caption{The comparison of different algorithms
        for inpainting incomplete RGB image. (From
        left to right) Observed images, images reconstructed by
        LRTI, images reconstructed by HaLRTC, images reconstructed by WTucker,
        images reconstructed by IRTD.
        Top row: $50\%$ missing.
        Middle row: $80\%$ missing.
        Bottom row: $90\%$ missing.}
    \label{fig:lena}
\end{figure*}

\subsection{Synthetic and Chemometrics Data}
In this subsection, we carry out experiments on synthetic and
chemometrics data. Two sets of synthetic data are generated and
both of them are $3$rd-order tensors of size $32\times 32\times
32$. The first tensor is generated according to the CP model which
is a summation of six equally-weighted rank-one tensors, with all
of the factor matrices drawn from a normal distribution. Thus the
truth rank is $6$ or $(6,6,6)$ in a multilinear rank form. The
other tensor is generated based on the Tucker decomposition model,
with a random core tensor of size $(3,4,5)$ multiplied by random
factor matrices along each mode. Clearly the groundtruth for the
multilinear rank is $(3,4,5)$. Two chemometrics data sets are also
considered in our simulations. One is the Amino Acid data set of
size $5\times 201\times 61$ and the other is the flow injection
data set of size $12\times 100\times 89$.




For each data set, we consider two cases with $50\%$ or $80\%$
entries missing in our simulations, where the missing entries are
randomly selected. The observed entries are corrupted by zero mean
Gaussian noise and the signal-to-noise ratio (SNR) is set to
$10$dB. The tensor reconstruction performance is evaluated by the
normalized mean squared error (NMSE)
$\|\boldsymbol{\mathcal{X}}-\boldsymbol{\mathcal{\hat{X}}}\|_F/
\|\boldsymbol{\mathcal{X}}\|_F$, where $\boldsymbol{\mathcal{X}}$
and $\boldsymbol{\mathcal{\hat{X}}}$ denote the true tensor and
the estimated one, respectively. The parameter $\lambda$ for the
LRTI is carefully selected as $\lambda=0.3$ for the synthetic data
set generated by the Tucker model and $\lambda=0.2$ for the other
data sets. The pre-defined multilinear rank for WTucker is set to
be $(12,12,12)$, $(6,8,10)$, $(5,10,10)$ and $(10,10,10)$ for the
CP, Tucker, Amino Acid and flow injection dataset, respectively.
For our proposed method, we choose $\lambda_1=0.5$ for all data
sets. Results are averaged over $10$ independent runs. The rank or
multilinear rank is estimated as the the most frequently occurring
rank or multilinear rank value. The standard deviation of the
estimated rank is also reported as an indication of the stability
of the inferred rank. Results are shown in Table \ref{table1}.
\begin{itemize}
    \item We observe that the proposed method presents the best recovery accuracy
    among all algorithms and can reliably estimate the true
    rank.
    \item Compared with the CP-decomposition based method LRTI, our proposed
    method
    presents a clear performance advantage over
    the LRTI when synthetic data are generated according to the Tucker model and
    slightly outperforms the LRTI even when data are generated according to the CP model.
    \item Our proposed method surpasses the tensor nuclear-norm based method HaLRTC by a big
    margin, which corroborates our claim that
    the proposed group log-sum functional is more effective than the tensor nuclear-norm
    in approximating the multilinear rank.
    \item Tucker model-based methods such as HaLRTC and WTucker work well
    only when synthetic data are generated according to
    the Tucker model, while our proposed method provides decent recovery
    performance whether data are
    generated via the Tucker or CP model.
    \item We also observe that
    our proposed method has similar run times as the other three algorithms.
    As the number of missing entries increases, our proposed method might need a few more iterations
    to reach convergence, and thus the average run time increases
    with the number of missing entries.
\end{itemize}
\subsection{Image Inpainting}
\begin{table}
    \centering
    \caption{The comparison of different algorithms for
        tensor completion (Image Inpainting)}
    \begin{tabular}{|c|c|c|c|}
        \hline
        & $50\%$ & $80\%$ & $90\%$ \\
        \hline
        LRTI & $0.0023$ & $0.0046$ & $0.0098$ \\
        \hline
        HaLRTC & $0.0019$ & $0.0066$ & $0.0134$ \\
        \hline
        WTucker & $0.0032$ & $0.0751$ & $0.2187$ \\
        \hline
        IRTD & $0.0015$ & $0.0046$ & $0.0082$ \\
        \hline
    \end{tabular}
    \label{table2}
\end{table}


The goal of image inpainting is to complete an image with missing
pixels. For a two-dimensional RGB picture, we can treat it as a
three-dimensional tensor. Here we consider imputing an incomplete
RGB image via respective algorithms. The benchmark Lena image is
used, with $50\%$, $80\%$ and $90\%$ missing entries considered in
our simulations. The recovery accuracy is evaluated by the MSE
metric which is defined as
$\text{MSE}=\frac{1}{M}\|\mathcal{O}^{\complement}\ast(\mathcal{X}-\hat{\mathcal{X}})\|_F^2$,
where $\mathcal{X}$ and $\hat{\mathcal{X}}$ respectively denote
the original normalized image and the estimated one,
$\mathcal{O}^{\complement}$ denotes the complement set of the
observed set, i.e. $\mathcal{O}$, and $M$ denotes the number of
missing elements. For LRTI, the parameter used to control the
tradeoff between the data fitting and rank is carefully selected
to $5$, $3$ and $3$ for $50\%$, $80\%$ and $90\%$ missing entries,
respectively.
The pre-defined rank for the WTucker is set to
$(100,100,3)$, $(80,80,3)$ and $(50,50,3)$ for $50\%$, $80\%$ and
$90\%$ missing entries, respectively. For our proposed method,
$\lambda_1$ is set to $3$, $0.5$ and $0.3$ for $50\%$, $80\%$ and
$90\%$ missing entries, respectively. The observed and recovered
images are shown in Fig.\ref{fig:lena} and MSEs of respective
algorithms are shown in Table \ref{table2}. From Table
\ref{table2}, we see that the proposed method renders a reliable
recovery even with $90\%$ missing entries, while the other two
Tuck model-based methods WTucker and HaLRTC may incur a
considerable performance degradation when the missing ratio is
high.


\section{Conclusions}\label{sec:conclusion}
In this paper, we proposed an iterative reweighted algorithm to
decompose an incomplete tensor into a concise Tucker
decomposition. To automatically determine the model complexity, we
introduced a new notion called order-$(N-1)$ sub-tensor and
introduced a group log-sum penalty on every order-$(N-1)$
sub-tensors to achieve a structural sparse core tensor. By
shrinking the zero order-$(N-1)$ sub-tensors, the core tensor
becomes a smaller one and a compact Tucker decomposition can be
obtained. By resorting to the majorization-minimization approach,
an iterative reweight algorithm was developed. Also, the
over-relaxed monotone fast iterative shrinkage-thresholding
technique is adapted and embedded in the iterative reweighted
process to reduce the computational complexity. The performance of
the proposed method is evaluated using synthetic data and real
data. Simulation results show that the proposed method offers
competitive performance compared with existing methods.

\bibliography{newbib}
\bibliographystyle{IEEEtran}

\end{document}